\documentclass[12pt]{amsart}
\usepackage[utf8]{inputenc}
\usepackage{pifont}
\usepackage[numbers]{natbib}
\usepackage{enumerate,verbatim,calc}
\usepackage{graphicx}
\usepackage[top=1.4in, bottom=1.4in, left=.9in, right=.9in]{geometry}
\usepackage[mathscr]{euscript}
\usepackage{mathrsfs}
\usepackage{amsmath,euscript,amsfonts,amsthm,amssymb,latexsym}

\def\R{{\mathbb R}}

\def\C{{\mathbb C}}
\def\D{{\mathbb D}}
\def\B{{\mathbb B}}

\def\L{{\mathbb L}}
\def\M{{\mathbb M}}

\def\h2{{\mathcal H}^2(\B^{2m})}
\def\b2{{\mathcal A}^2(\B^{2m})}

\begin{document}
\title{\bf A note on Cartan isometries}

\author{Ameer Athavale}
\address{Department of Mathematics, Indian Institute of Technology Bombay, Powai, Mumbai 400076, India}
\email{athavale@math.iitb.ac.in}
\begin{abstract}
We record a lifting theorem for the intertwiner of two $S_{\Omega}$-isometries which are those subnormal operator tuples whose minimal normal extensions have their Taylor spectra contained in the Shilov boundary of a certain function algebra associated with $\Omega$, $\Omega$ being a bounded convex domain in $\C^n$ containing the origin. The theorem captures several known lifting results in the literature and yields interesting new examples of liftings as a consequence of its being applicabile to Cartesian products $\Omega$ of classical Cartan domains in $\C^n$. Further, we derive intrinsic characterizations of $S_{\Omega}$-isometries where $\Omega$ is a classical Cartan domain of any of the types I, II, III and IV, and we also provide a neat description of an $S_{\Omega}$-isometry in case $\Omega$ is a finite Cartesian product of such Cartan domains.
\end{abstract}
\subjclass[2010]{Primary 47A13, 47B20}
\keywords{Cartan domain, Cartan isometry, spherical isometry}
\maketitle
\section{Introduction}
For $\mathcal H$ a complex infinite-dimensional separable Hilbert space, we use ${\mathcal B}({\mathcal H})$ to denote the algebra of bounded linear operators on $\mathcal H$. 
An $n$-tuple $S = (S_1,\ldots, S_n)$ of commuting operators $S_i$ in ${\mathcal B}({\mathcal H})$ is said to be {\it subnormal} if there exist a Hilbert
space ${\mathcal K}$ containing ${\mathcal H}$ and an $n$-tuple $N = (N_1, \ldots,
N_n)$ of commuting normal operators $N_i$ in ${\mathcal B}({\mathcal K})$ such
that $N_i {\mathcal H}  \subset {\mathcal H}$ and $N_i/{\mathcal H} = S_i$ for $1
\leq i \leq n$.  \\

Suppose $S = (S_1,\ldots, S_n)$ is a tuple of commuting operators in 
${\mathcal B}({\mathcal H})$ and
$T = (T_1,\ldots, T_n)$ a tuple of commuting operators in ${\mathcal B}({\mathcal J})$. If there exists a bounded linear operator $X : {\mathcal H} \rightarrow {\mathcal J}$ such that $XS_i = T_iX$ for each $i$, then $X$ is said to be an {\it intertwiner} (for $S$ and $T$) and we denote this fact by $XS = TX$. If $X : {\mathcal H} \rightarrow {\mathcal J}$ and
$Y : {\mathcal J} \rightarrow {\mathcal H}$ are two intertwiners for $S$ and $T$ such that $XS = TX$ and $YT=SY$, and both $X$ and $Y$ are injective and have dense ranges, then $S$ is said to be {\it quasisimilar to} $T$. The operator tuple $S$ is said to be  {\it unitarily equivalent to $T$} if one can find a unitary intertwiner for $S$ and $T$. Any subnormal operator tuple is known to admit a `minimal' normal extension that is unique up to unitary equivalence (see \cite {I}).\\

For a bounded domain $\Omega$ in $\C^n$, we let $A(\Omega) = \{f \in C({\bar {\Omega}}): f {\rm \ is \ holomorphic \ on\ } \Omega\}$, where $C({\bar {\Omega}})$ denotes the algebra of continuous functions on the closure ${\bar {\Omega}}$ of $\Omega$. The {\it Shilov boundary} of $A(\Omega)$ (or $\Omega$) is defined to be the smallest closed subset $S_{\Omega}$ of 
${\bar {\Omega}}$ such that, for any $f \in A(\Omega)$,
$$ \sup\{|f(z)|: z \in {\bar {\Omega}}\} = \sup\{|f(z)|: z \in S_{\Omega}\}.$$\\

Of special interest to us are domains $\Omega$ that are Cartesian products $\Omega_1\times\cdots\times \Omega_m$ with $\Omega_i \subset \C^{n_i}$ being a classical Cartan domain of any of the four types I II, III and IV (refer to \cite{Car}, \cite{Hu},  \cite{J-P}, \cite{L}); any such domain $\Omega$ will be referred to as a {\it standard Cartan domain}. The open unit ball $\B_n$ in $\C^n$ is a classical Cartan domain of type I with its Shilov boundary coinciding with the unit sphere in $\C^n$. The open unit polydisk $\D^n$ in $\C^n$ is a standard Cartan domain with its Shilov boundary coinciding with the unit polycircle in $\C^n$. The standard Cartan domains are special examples of bounded symmetric domains and are `circled around the origin' in the sense that they contain the origin and are invariant under multiplication by $e^{\sqrt{-1}\theta}$, $\theta \in \R$. It follows from \cite[Lemma 5.7]{Di-E} that the Shilov boundary $S_{\Omega}$ of any standard Cartan domain $\Omega = \Omega_1\times\cdots  \times \Omega_m$, where each $\Omega_i$ is a classical Cartan domain in $\C^{n_i}$, is given by $S_{\Omega}= S_{\Omega_1}\times\cdots\times S_{\Omega_m}$.  \\ 

A subnormal tuple $S$ will be referred to as an {\it $S_{\Omega}$-isometry} if the Taylor spectrum $\sigma(N)$ of its minimal normal extension $N$ is contained in the Shilov boundary $S_{\Omega}$ of $\Omega$. We use $I_{\mathcal H}$ (resp. $0_{\mathcal H}$) to denote the identity operator (resp. the zero operator) on $\mathcal H$. An $S_{\B_n}$-isometry is precisely a {\it spherical isometry}, that is, an $n$-tuple $S$ of commuting operators $S_i$ in
${\mathcal B}({\mathcal H})$ satisfying $\sum_{i=1}^nS_i^*S_i = I_{\mathcal H}$ (refer to \cite [Proposition 2]{At1}). An $S_{\D^n}$-isometry is precisely a {\it toral isometry}, that is, an $n$-tuple $S$ of commuting operators $S_i$ in ${\mathcal B}({\mathcal H})$ satisfying $S_i^*S_i = I_{\mathcal H}$ for each $i$ (refer to \cite [Proposition 6.2]{N-F}).  Any $S_{\Omega}$-isomerty with $\Omega$ a standard Cartan domain will be referred to as a {\it Cartan isometry}.\\

We will say that a domain $\Omega \subset \C^n$ {\it satisfies the property (A)} if, for any positive regular Borel measure $\eta$ supported on the Shilov boundary $S_{\Omega}$ of $\Omega$, the triple $(A(\Omega)|S_{\Omega},S_{\Omega},\eta)$ is {\it regular} in the sense of \cite{A}, that is, for any positive continuous function $\phi$ defined on $S_{\Omega}$, there exists a sequence of functions $\{\phi_m\}_{m \geq 1}$ in $A(\Omega)$ such that $|\phi_m| < \phi$ on $S_{\Omega}$ and $\lim_{m \rightarrow \infty} |\phi_m| = \phi$ $\eta$-almost everywhere. \\

The discussion in Section 5 of \cite {Di-E} shows that any bounded symmetric domain circled around the origin satisfies the property (A).\\

In Section 2, we state a lifting result for the intertwiner of certain $S_{\Omega}$-isometries of which Cartan isometries are special examples. In Section 3 we provide an intrinsic characterization of $S_{\Omega}$-isometries for Cartan domains $\Omega$ of type IV and then characterize $S_{\Omega}$-isometries for $\Omega$ a Cartesian product of the open unit balls and Cartan domains of type IV (see Theorem 3.5). In Section 4, we characterize $S_{\Omega}$-isometries for Cartan domains of type I and observe that Theorem 3.5 holds with the open unit balls replaced by Cartan domains of type I. Finally, in Section 5 we characterize $S_{\Omega}$-isometries for Cartan domains of type II and of type III and end up with a substantial generalization of Theorem 3.5. For basic facts pertaining to classical Cartan domains and bounded symmetric domains in general, the reader is referred to \cite{Hu}, \cite{J-P} and \cite{L}. It may be noted that Shilov boundaries are referred to as `characteristic manifolds' in \cite{Hu}.

\section{A lifting theorem for certain $S_{\Omega}$-isometries}
The proof of Theorem 2.1 below is similar to the proofs of \cite[Theorem 3.2]{At2} and \cite[Proposition 4.6]{At3}; however, unlike there, it circumvents using the Taylor functional calculus of \cite{T}. Also, unlike in \cite{At2} and \cite{At3}, the Shilov boundary $S_{\Omega}$ of $\Omega$ may not coincide with the topological boundary $\partial \Omega$ of $\Omega$.\\

{\bf Theorem 2.1}. Let $\Omega$ be a bounded convex domain in $\C^n$ containing the origin and satisfying the property (A) of Section 1. Let $S =(S_1,\ldots,S_n) \in {\mathcal B}({\mathcal H})^n$ and $T =(T_1,\ldots,T_n) \in {\mathcal B}({\mathcal J})^n$ be $S_{\Omega}$-isometries, and let 
$M =(M_1,\ldots,M_n) \in {\mathcal B}({\tilde {\mathcal H}})^n$ and $N =(N_1,\ldots,N_n) \in {\mathcal B}({\tilde {\mathcal J}})^n$ respectively be the minimal normal extensions of $S$ and $T$. If $X: {\mathcal H} \rightarrow {\mathcal J}$ is an intertwiner for $S$ and $T$, then $X$ lifts to a (unique) intertwiner $\tilde X: {\tilde {\mathcal H}} \rightarrow {\tilde {\mathcal J}}$ for $M$ and $N$; moreover, $\|{\tilde X}\| =\|X\|$.
\begin{proof}
Let $f \in A(\Omega)$. For any positive integer $m \geq 2$, $f_m$ defined by $f_m(z) = f((1-\frac{1}{m})z)$ is holomorphic on an open neighborhood of $\bar \Omega$. 
Since $\bar \Omega$ is polynomially convex, $f_m$ is the uniform limit (on $\bar\Omega$) of a sequence $\{p_{m,k}\}_k$ of polynomials by the Oka-Weil approximation theorem (see \cite{Ra}, Chapter VI, Theorem 1.5). If $X$ intertwines $S$ and $T$, then one clearly has $Xp_{m,k}(S) = p_{m,k}(T)X$. 
If $\rho_M$ and $\rho_N$ are respectively the spectral measures of $M$ and $N$ (supported on $S_{\Omega}$), then $\rho_S= P_{\mathcal H}\rho_M|{\mathcal H}$ and $\rho_T= P_{\mathcal J}\rho_N|{\mathcal J}$ are respectively the semi-spectral measure of $S$ and $T$ with $P_{\mathcal H}$ and $P_{\mathcal J}$ being appropriate projections, and for any $u \in \mathcal H$ and any $v \in \mathcal K$ one has
$$ \|p_{m,k}(S)u\|^2 = \int_{S_{\Omega}}|p_{m,k}(z)|^2 d\langle\rho_S(z)u,u\rangle$$
and
$$ \|p_{m,k}(T)v\|^2 = \int_{S_{\Omega}}|p_{m,k}(z)|^2 d\langle\rho_T(z)v,v\rangle.$$
Choosing $v = Xu$ and using $Xp_{m,k}(S) = p_{m,k}(T)X$, one has
$$ \int_{S_{\Omega}}|p_{m,k}(z)|^2 d\langle\rho_T(z)Xu,Xu\rangle \leq \|X\|^2\int_{S_{\Omega}}|p_{m,k}(z)|^2 d\langle\rho_S(z)u,u\rangle.$$
Letting first $k$ tend to infinity and then $m$ tend to infinity, one obtains
$$ \int_{S_{\Omega}}|f(z)|^2 d\langle\rho_T(z)Xu,Xu\rangle \leq \|X\|^2\int_{S_{\Omega}}|f(z)|^2 d\langle\rho_S(z)u,u\rangle.$$
Consider $\eta(\cdot) = \langle\rho_T(\cdot)Xu,Xu\rangle + \langle\rho_S(\cdot)u,u\rangle$. Since $(A(\Omega)|S_{\Omega},S_{\Omega},\eta)$ is a regular triple, for any positive continuous function $\phi$ on $S_{\Omega}$ there exists a sequence of functions $\{\phi_m\}_{m \geq 1}$ in $A(\Omega)$ such that $|\phi_m| < \sqrt \phi$ on $S_{\Omega}$ and $\lim_{m \rightarrow \infty} |\phi_m| = \sqrt \phi$ $\eta$-almost everywhere. Replacing $f$ by $\phi_m$ in the last integral inequality and letting $m$ tend to infinity, one obtains
$$ \int_{S_{\Omega}}\phi(z) d\langle\rho_T(z)Xu,Xu\rangle \leq \|X\|^2\int_{S_{\Omega}}\phi(z) d\langle\rho_S(z)u,u\rangle.$$
That yields $\langle\rho_T(\cdot)Xu,Xu\rangle \leq \|X\|^2\langle\rho_S(\cdot)u,u\rangle$ for every $u$ in $\mathcal H$. The desired conclusion now follows by appealing to \cite[Lemma 4.1]{M}. 
\end{proof}

In so far as the function algebra $A(\Omega)$ is concerned, Theorem 2.1 is an improvement over \cite[Theorem 5.1]{M} by virtue of its using the more widely applicable property (A) in place of the property `approximating in modulus' as required of a function algebra in \cite{M}.\\ 

{\bf Corollary 2.2}. Let $\Omega$ be any bounded symmetric domain circled around the origin (so that $\Omega$ can in particular be a standard Cartan domain). Let $S =(S_1,\ldots,S_n) \in {\mathcal B}({\mathcal H})^n$ and $T =(T_1,\ldots,T_n) \in {\mathcal B}({\mathcal J})^n$ be $S_{\Omega}$-isometries, and let 
$M =(M_1,\ldots,M_n) \in {\mathcal B}({\tilde {\mathcal H}})^n$ and $N =(N_1,\ldots,N_n) \in {\mathcal B}({\tilde {\mathcal J}})^n$ respectively be the minimal normal extensions of $S$ and $T$. If $X: {\mathcal H} \rightarrow {\mathcal J}$ is an intertwiner for $S$ and $T$, then $X$ lifts to a (unique) intertwiner $\tilde X: {\tilde {\mathcal H}} \rightarrow {\tilde {\mathcal J}}$ for $M$ and $N$; moreover, $\|{\tilde X}\| =\|X\|$.
\begin{proof}
Any bounded symmetric domain circled around the origin is convex by \cite[Corollary 4.6]{L} and, as noted in Section 1, satisfies the property (A).
\end{proof}

{\bf Remark 2.3}. Letting $\Omega$ to be the open unit ball $\B_n$ in $\C^n$, Corollary 2.2 captures \cite [Proposition 8]{At1} which is a lifting result for the intertwiner of spherical isometries. Letting $\Omega$ to be the open unit polydisk $\D^n$ in $\C^n$, Corollary 2.2 captures  \cite[Proposition 5.2]{M} which is a lifting result for the intertwiner of toral isometries.  In \cite{At3}, the author introduced a class $\Omega^{(n)}$ of convex domains $\Omega_p$ in $\C^n$ that satisfy the property (A); for $n \geq 2$, the class $\Omega^{(n)}$ happens to be distinct from the class of strictly pseudoconvex domains and the class of bounded symmetric domains in $\C^n$. Letting $\Omega$ to be $\Omega_p$, Theorem 2.1 (but not Corollary 2.2) captures \cite[Proposition 4.6]{At3}. A variant of Theorem 2.1 that is valid for (not necessarily convex) strictly pseudoconvex bounded domains $\Omega$ with $C^2$ boundary was proved in \cite{At2}; however, Theorem 2.1 does apply to strictly pseudoconvex bounded domains that are convex since any strictly pseudoconvex bounded domain $\Omega$ is known to satisfy the property (A) (refer to \cite{A} and \cite{Di-E}).
\\

{\bf Remark 2.4}. Arguing as in \cite[Theorem 5.2]{M}, one can establish the following facts in the context of Theorem 2.1: If $X$ is isometric, then so is ${\tilde X}$; if $X$ has dense range, then so has ${\tilde X}$; if $X$ is bijective, then so is ${\tilde X}$.  Also, it follows from \cite[Lemma 1]{At1} that if $S$ and $T$ of Theorem 2.1 are quasisimilar, then the minimal normal extensions of $S$ and $T$ are unitarily equivalent (cf. \cite[Proposition 9]{At1}).\\

\section{Lie sphere isometries: $S_{\Omega}$-isometries for Cartan domains $\Omega$ of type IV}
The Lie ball $\L_n$ in $\C^n$ is defined by
$$ \L_n =\left \{z \in \C^n: \left(\|z\|^2+\sqrt{\|z\|^4-|\langle z,{\bar z}\rangle|^2}\right)^{1/2} < 1 \right\}.$$

Lie balls are classical Cartan domains $\Omega_{IV}(n)$. We note that $\L_1 = \D^1 = \B_1$. The Shilov boundary $S_{\L_n}$ of $\L_n$ (also referred to as the {\it Lie sphere}) is given by 
$$ S_{\L_n} =\{(z_1,\ldots,z_n): z_i =x_ie^{\sqrt{-1}\theta}(1 \leq i \leq n),\ \theta \in \R,\  x_i \in \R,\  x_1^2+\cdots +x_n^2=1\}.$$

We will refer to an $S_{\L_n}$-isometry as a {\it Lie sphere isometry}; thus Lie sphere isometries are $S_{\Omega}$-isometries for classical Cartan domains $\Omega$  of type IV. It should be noted that $S_{\L_n}$ is contained in $S_{\B_n}$ so that any Lie sphere isometry is a spherical isometry! We plan to provide an intrinsic characterization of a Lie sphere isometry, and for that purpose we need Lemma 3.1 below. (A result more general than that of Lemma 3.1 is present in the unpublished work \cite{Ch}; we present here a direct proof for the reader's convenience).\\

{\bf Lemma 3.1}. Let $S =(S_1,\ldots,S_n) \in {\mathcal B}({\mathcal H})^n$ be a subnormal tuple with the minimal normal extension $N =(N_1,\ldots,N_n) \in {\mathcal B}({\mathcal K})^n$. If $S_i^*S_j =S_j^*S_i$ (so that  $S_i^*S_j$ is self-adjoint) for some $i$ and $j$, then $N_i^*N_j =N_j^*N_i$ (so that $N_i^*N_j$ is also self-adjoint).
\begin{proof}
For arbitrary non-negative integers $k_i$ and $l_i$ ($1 \leq i \leq n$), consider 
$$\langle (N_i^*N_j -N_j^*N_i)({N_1^*}^{k_1}\cdots {N_n^*}^{k_n}x),({N_1^*}^{l_1}\cdots{ N_n^*}^{l_n}y)\rangle \ \ (x,y \in {\mathcal H}).$$
Using that $N_p$ and $N_q^*$ commute for all $p$ and $q$ and $N_p|{\mathcal H} =S_p$ for every $p$, it is easy to see that this inner product reduces to
$$\langle (S_i^*S_j -S_j^*S_i)({S_1}^{l_1}\cdots {S_n}^{l_n}x),({S_1}^{k_1}\cdots{S_n}^{k_n}y)\rangle.$$ 
Since ${\mathcal K}$ is the closed linear span of vectors of the type ${N_1^*}^{k_1}\cdots {N_n^*}^{k_n}x$, the desired result is obvious.
\end{proof}

{\bf Theorem 3.2}. For an $n$-tuple $S =(S_1,\ldots,S_n)$ of operators $S_i$ in ${\mathcal B}({\mathcal H})$,  (a) and (b) below are equivalent.\\
(a) $S$ is a Lie sphere isometry.\\
(b) $S$ is a spherical isometry and $S_i^*S_j$ is self-adjoint for every $i$ and $j$.
\begin{proof}
Suppose (a) holds so that $S =(S_1,\ldots,S_n) \in {\mathcal B}({\mathcal H})^n$ is a Lie sphere isometry. Then the minimal normal extension $N =(N_1,\ldots,N_n) \in {\mathcal B}({\mathcal K})^n$ of $S$ has its Taylor spectrum $\sigma(N)$ contained in $S_{\L_n}$. Since for any $(z_1,\ldots,z_n) \in S_{\L_n}$ the equalities $|z_1|^2+\cdots +|z_n|^2=1$ and ${\bar z_i}z_j-{\bar z_j}z_i=0\ (1 \leq i,j \leq n)$ hold, one has $N_1^*N_1+\cdots+N_n^*N_n = I_{\mathcal K}$ and $N_i^*N_j -N_j^*N_i =0_{\mathcal K} \ (1 \leq i,j \leq n)$. Compressing these  equations to $\mathcal H$, (b) is seen to hold.

Conversely, suppose (b) holds. Since one has $\sum_iS_i^*S_i =I_{\mathcal H}$, \cite[Proposition 2]{At2} gives that $S$ is a subnormal tuple with the Taylor spectrum $\sigma(N)$ of its minimal normal extension $N$ contained in the unit sphere $S_{\B_n}$. The condition that $S_i^*S_j$ is self-adjoint for every $i$ and $j$ guarantees, by Lemma 3.1, that $N_i^*N_j -N_j^*N_i =0_{\mathcal K}$ for every $i$ and $j$. It follows then from the spectral theory for $N$ that the Taylor spectrum of $N$ is contained in the set $\{z\in S_{\B_n}: {\bar z_i}z_j-{\bar z_j}z_i=0 {\rm \  for \ every \ } i {\rm\  and \ } j\}$ which, as an elementary verification using polar coordinates shows, is the set $S_{\L_n}$.
\end{proof}

At this stage we introduce a notational convention that will be convenient to use in the sequel. For a complex polynomial $p(z,w) = \sum_{\alpha, \beta}a_{\alpha,\beta}z^{\alpha}w^{\beta}$ in the variables $z,w \in \C^n$ and for any $n$-tuple $S$ of commuting operators $S_i$ in ${\mathcal B}({\mathcal H})$, $p(z,w)(S,S^*)$ is to be interpreted as $\sum_{\alpha,\beta}a_{\alpha,\beta}{S^*}^{\beta}S^{\alpha}$. Thus $S$ is a spherical isometry if and only if $(1-\sum_{i=1}^nz_iw_i)(S,S^*)=0_{\mathcal H}$. A {\it contraction} is an operator $S$ in ${\mathcal B}({\mathcal H})$ for which $(I-S^*S) \equiv (1-zw)(S,S^*) \geq 0_{\mathcal H}$. As proved in \cite{At0}, an $n$-tuple $S$ of commuting contractions $S_i$ in ${\mathcal B}({\mathcal H})$ is subnormal if and only if  $\Pi_{i=1}^n (1-z_iw_i)^{k_i}(S,S^*) \geq 0_{\mathcal H}$ for all non-negative integers $k_i$. Further, with $p(z,w)$ as here and with $S$ a subnormal tuple, the proof of Lemma 3.1 goes through with $S_i^*S_j-S_j^*S_i$ there replaced by $p(z,w)(S,S^*)$. We state this generalization (due to Chavan) of \cite[Proposition 8] {At-P} as Lemma 3.3.\\

{\bf Lemma 3.3 \cite{Ch}}. Let $S \in {\mathcal B}({\mathcal H})^n$ be a subnormal tuple with the minimal normal extension $N \in {\mathcal B}({\mathcal K})^n$. If $p(z,w)(S,S^*) =0_{\mathcal H}$, then $p(z,w)(N,N^*)=0_{\mathcal K}$.\\

{\bf Lemma 3.4}. Let $S=(S_1,\ldots,S_n)$ be a tuple of commuting operators in ${\mathcal B}({\mathcal H})$ such that each $S_i$ is a coordinate of a subtuple of $S$ that is a spherical isometry. Then $S$ is subnormal.
\begin{proof} Suppose for each $i$ there exist positive integers $j(i,1),\ldots, j(i,p_i)$, with $j(i,k) = i$ for some $k$, such that $(S_{j(i,1)},\ldots, S_{j(i,k)}=S_i, \ldots, S_{j(i,p_i)})$ is a spherical isometry. It is clear that each $S_i$ is then a contraction. We need to verify that $\Pi_{i=1}^n (1-z_iw_i)^{k_i}(S,S^*) \geq 0$ for all non-negative integers $k_i$. The verification results by writing
each factor $(1-z_iw_i)$ as
$$ (1-z_iw_i) = (\{1-\sum_{l=1}^{p_i}z_{j(i,l)}w_{j(i,l)}\} + \sum_{\substack{l=1\\ l \neq k}}^{p_i}z_{j(i,l)}w_{j(i,l)}).$$
\end{proof}

We are now in a position to characterize $S_{\Omega}$-isometries in case $\Omega$ is a Cartesian product of the open unit balls and the Lie balls. A substantial generalization of Theorem 3.5 below will be achieved in Section 5; however, the essential ingredients of the relevant argument are present in the proof of Theorem 3.5 and occur at this stage without the clutter of too many ideas.\\

{\bf Theorem 3.5}. Let $\Omega = \Omega_1\times\cdots\times\Omega_m \subset \C^{n}$ where each $\Omega_i$ is either the open unit ball in $\C^{n_i}$ or the Lie ball in $\C^{n_i}$ (and where $n =n_1+\cdots+n_m$). Let $S_i=(S_{i,1},\ldots,S_{i,n_i})$ be an $n_i$-tuple of operators in ${\mathcal B}({\mathcal H})$ for $1 \leq i \leq m$ and let the operator coordinates of the $n$-tuple $S=(S_1;\ldots;S_m)$ commute with each other. Then $S$ is an $S_{\Omega}$-isometry if and only if each $S_i$ is an $S_{\Omega_i}$-isometry.
\begin{proof}
We illustrate the proof for the case $m=2$, $n_1 =2$, $n_2=3$, $\Omega_1 =\B{_2}$ and $\Omega_2 = \L{_3}$. The general case is then no more than an exercise in notational book-keeping.

Suppose first that $S=(S_1;S_2)=(S_{1,1},S_{1,2};S_{2,1},S_{2,2},S_{2,3})$ is an $S_{\B_2\times\L_3}$-isometry so that $S$ is subnormal and the Taylor spectrum $\sigma(N)$ of its minimal normal extension $N=(N_1;N_2)=(N_{1,1},N_{1,2};N_{2,1},N_{2,2},N_{2,3}) \in {\mathcal B}({\mathcal K})^5$ is contained in $S_{\B_2\times\L_3}= S_{\B_2}\times S_{\L_3}$. By the projection property of the Taylor spectrum (refer to \cite{T}), the inclusions $\sigma(N_1) \subset S_{\B_2}$ and $\sigma(N_2)\subset S_{\L_3}$ hold. While $N_1$ and $N_2$ may not be the minimal normal extensions of $S_1$ and $S_2$, they certainly satisfy the relations
$$ \sum_{i=1}^2N_{1,i}^*N_{1,i} = I_{\mathcal K}, \sum_{j=1}^3N_{2,j}^*N_{2,j} = I_{\mathcal K}, N_{2,k}^*N_{2,l}=N_{2,l}^*N_{2,k}, 1 \leq k,l \leq 3.$$ 
Compressing these equations to $\mathcal H$, one obtains
$$ \sum_{i=1}^2S_{1,i}^*S_{1,i} = I_{\mathcal H}, \sum_{j=1}^3S_{2,j}^*S_{2,j} = I_{\mathcal H}, S_{2,k}^*S_{2,l}=S_{2,l}^*S_{2,k}, 1 \leq k,l \leq 3.$$ 
Using our observations in Section 1 related to spherical isometries and appealing to Theorem 3.2, it  follows that $S_1$ is  an $S_{\B_2}$-isometry and $S_2$ is an  $S_{\L_3}$-isometry.

Conversely, suppose $S_1=(S_{1,1},S_{1,2})$ is an $S_{\B_2}$-isometry and $S_2= (S_{2,1},S_{2,2},S_{2,3})$ is an $S_{\L_3}$-isometry. Then the identities for $S$ as recorded above hold so that
$$ (1-\sum_{i=1}^2z_iw_i)(S_1,{S_1}^*) =0_{\mathcal H}, (1-\sum_{j=1}^3z_jw_j)(S_2,{S_2}^*) =0_{\mathcal H}, (z_lw_k-z_kw_l)(S_2,{S_2}^*)=0_{\mathcal H}, 1 \leq k,l \leq 3.$$ 
While both $S_1$ and $S_2$ are subnormal, the crucial thing to verify  is that $S=(S_1;S_2)$ is subnormal. But the subnormality of $S$ is now a consequence  of Lemma 3.4.
Letting $N= (N_1;N_2)$ to be the minimal normal extension of $S=(S_1;S_2)$ and using Lemma 3.3, we see that $N$  satisfies the same identities as $S$. That $\sigma(N)$ is contained in  $S_{\B_2}\times S_{\L_3} = S_{\B_2 \times \L_3}$ is now a consequence of the spectral theory for $N$. 
\end{proof}
\section{$S_{\Omega}$-isometries for Cartan domains $\Omega$ of type I}
We use the symbol $\M(p,q)$ to denote the set of  complex matrices of order $p\times q$ and the symbol $I_n$ to denote the identity matrix of order $n$. The complex tranjugate of a matrix $Z$ will be denoted by $Z^*$ so that $Z^*$ is the transpose ${\bar Z}^t$ of the complex conjugate $\bar Z$ of $Z$. The classical Cartan domain $\Omega_{I}(p,q)$ of type I in $\C^n$ is defined by the following conditions:
$$n=pq,\ 1 \leq p \leq q, \ \Omega_{I}(p,q) = \{Z \in \M(p,q): I_p -ZZ^* \geq 0\}$$
The Shilov boundary of  $\Omega_{I}(p,q)$ is given by
$$ S_{\Omega_{I}(p,q)} = \{Z \in \M(p,q): I_p -ZZ^* = 0\}.$$
It will be convenient to rewrite $\Omega_{I}(p,q)$ as
$$ \Omega_{I}(p,q) = \{(z_{1,1},\ldots, z_{1,q};z_{2,1},\ldots, z_{2,q};\ldots;z_{p,1},\ldots, z_{p,q})\in \C^{pq}:
I_p -(z_{i,j})(\overline{z_{j,i}}) \geq 0\}$$
and $S_{\Omega_{I}(p,q)}$ as
$$ S_{\Omega_{I}(p,q)}= \{(z_{1,1},\ldots, z_{1,q};z_{2,1},\ldots, z_{2,q};\ldots;z_{p,1},\ldots, z_{p,q})\in \C^{pq}:
I_p -(z_{i,j})(\overline{z_{j,i}}) = 0\}.$$
The conditions defining the Shilov boundary $S_{\Omega_{I}(p,q)}$ can be written as
$$\sum_{k=1}^q { \overline{z_{j,k}} }z_{i,k} = \delta_{i,j},\ \ 1 \leq i \leq j \leq p.$$
Formally replacing $z_{i,j}$ by $S_{i,j}$ and ${\overline {z_{i,j}}}$ by $S_{i,j}^*$ (where $S_{i,j} \in {\mathcal B}({\mathcal H}))$, one is led to
$$\sum_{k=1}^q S_{j,k}^*S_{i,k} = \delta_{i,j}I_{\mathcal H},\ \ 1 \leq i \leq j \leq p.$$
{\bf Theorem 4.1}. For $p \leq q$, let $S_i=(S_{i,1},\ldots,S_{i,q})$ be a $q$-tuple of operators in ${\mathcal B}({\mathcal H})$ for $1 \leq i \leq p$  and let the operator coordinates of the $pq$-tuple $S=(S_1;\ldots;S_p)$ commute with each other. Then (a) and (b) below are equivalent.\\
(a) $S$ is an $S_{\Omega_{I}(p,q)}$-isometry.\\
(b) $$\sum_{k=1}^q S_{j,k}^*S_{i,k} = \delta_{i,j}I_{\mathcal H},\ 1 \leq i \leq j \leq p.$$
\begin{proof}
Suppose $S$ is an $S_{\Omega_{I}(p,q)}$-isometry. Then its minimal normal extension $N=(N_1;\ldots;N_p) \in {\mathcal B}({\mathcal K})^{pq}$ (with $N_i = (N_{i,1},\ldots,N_{i,q})$ for each $i$) has its Taylor spectrum $\sigma(N)$ contained in $S_{\Omega_{I}(p,q)}$. Since for any $z =(z_{1,1},\ldots,z_{p,q}) \in S_{\Omega_{I}(p,q)}$ the equalities $\sum_{k=1}^q {\overline {z_{j,k}}}z_{i,k} = \delta_{i,j},\  1 \leq i \leq j \leq p$ hold, one has $\sum_{k=1}^q N_{j,k}^*N_{i,k} = \delta_{i,j}I_{\mathcal K},\ \ 1 \leq i \leq j \leq p$. Compressing the last equations to $\mathcal H$, (b) is seen to hold.

Conversely, suppose (b) holds. The conditions in (b) corresponding to $1 \leq i=j \leq p$ guarantee that each $S_i$ is a spherical isometry. It then follows from Lemma 3.4 that $S=(S_1;\ldots;S_p)$ is subnormal. If $N$ in the notation used above is the minimal normal extension of $S$, then Lemma 3.3 yields the equalities 
$\sum_{k=1}^q N_{j,k}^*N_{i,k} = \delta_{i,j}I_{\mathcal K},\  1 \leq i \leq j \leq p$. The spectral theory for $N$ now implies that $\sigma(N)$ is contained in $S_{\Omega_{I}(p,q)}$.
\end{proof}

Using Theorems 3.2 and 4.1 and arguing as in the proof of Theorem 3.5, one can now establish Theorem 4.2 below.\\

{\bf Theorem 4.2}. Let $\Omega = \Omega_1\times\cdots\times\Omega_m \subset \C^{n}$ where each $\Omega_i$ is a classical Cartan domain of any of the types I and IV in $\C^{n_i}$ (and where $n =n_1+\cdots+n_m$). Let $S_i=(S_{i,1},\ldots,S_{i,n_i})$ be an $n_i$-tuple of operators in ${\mathcal B}({\mathcal H})$ for $1 \leq i \leq m$ and let the operator coordinates of the $n$-tuple $S=(S_1;\ldots;S_m)$ commute with each other. Then $S$ is an $S_{\Omega}$-isometry if and only if each $S_i$ is an $S_{\Omega_i}$-isometry.\\

{\bf Remark 4.3}. Since $\Omega_{1,n}$ is the open unit ball in $\C^n$, Theorem 4.1 generalizes the well-known characterization of an $S_{\B_n}$-isometry as a spherical isometry, the case $n=1$ of course yielding  the identification of an $S_{\B_1}$-isometry with an isometry. Also, Theorem 4.2 generalizes Theorem 3.5 and, with $\Omega_i$ chosen to be the unit disk $\D^1=\B_1$ in $\C$ for each $i$, yields the well-known characterization of an $S_{\D^n}$-isometry as a toral isometry. 

\section{$S_{\Omega}$-isometries for Cartan domains $\Omega$ of type II and of type III}
Let ${\mathcal S}(p) = \{Z \in \M(p,p): Z^t = Z\}$ and let ${\mathcal A}(p) = \{Z \in \M(p,p): Z^t = -Z\}$. The classical Cartan domain $\Omega_{II}(p)$ of type II in $\C^n$ is defined by the following conditions:
$$n=p(p+1)/2,\ p \geq 1,\ \Omega_{II}(p) = \{Z \in {\mathcal S}(p): I_p -ZZ^* \geq 0\}$$
The classical Cartan domain $\Omega_{III}(p)$ of type III in $\C^n$ is defined by the following conditions:
$$n=p(p-1)/2,\ p \geq 2,\ \Omega_{III}(p) = \{Z \in {\mathcal A}(p): I_p -ZZ^* \geq 0\}$$
(Some authors may refer to type II domains as type III domains and vice versa).\\
 
The Shilov boundary of  $\Omega_{II}(p)$ is given by
$$ S_{\Omega_{II}(p)} = \{Z \in {\mathcal S}(p): I_p -ZZ^* = 0\}$$
and the Shilov boundary of  $\Omega_{III}(2p)$ is given by
$$ S_{\Omega_{III}(2p)} = \{Z \in {\mathcal A}(2p): I_{2p} -ZZ^* = 0\}.$$
(We will comment on $S_{\Omega_{III}(2p+1)}$ later.)\\

We let  $$z_{{\mathcal S}(p)} = (z_{1,1},\ldots,z_{1,p};z_{2,2},\ldots,z_{2,p};\ldots;z_{p,p})$$
and 
$$z_{{\mathcal A}(p)} = (z_{1,2},\ldots,z_{1,p};z_{2,3},\ldots,z_{2,p};\ldots;z_{p-1,p}).$$
It will be convenient to rewrite $\Omega_{II}(p)$ as
$$ \Omega_{II}(p) =\{z_{{\mathcal S}(p)}\in \C^{p(p+1)/2}: {\rm With\ } z_{j,i} := z_{i,j} {\rm \ for\ } i \leq j, I_p -(z_{i,j})(\overline{z_{j,i}}) \geq 0\}$$
and $\Omega_{III}(p)$ as
$$ \Omega_{III}(p) =\{z_{{\mathcal A}(p)} \in \C^{p(p-1)/2}: {\rm With\ } z_{j,i} := -z_{i,j} {\rm \ for\ } i \leq j, I_p -(z_{i,j})(\overline{z_{j,i}}) \geq 0\}.$$
The conditions defining the Shilov boundary $S_{\Omega_{II}(p)}$ can  be written as follows:
$${\rm With\ } z_{j,i} := z_{i,j} {\rm \ for\ } i \leq j, \ \sum_{k=1}^p { \overline{z_{j,k}} }z_{i,k} = \delta_{i,j},\ \ 1 \leq i \leq j \leq p$$
Also, the conditions defining the Shilov boundary $S_{\Omega_{III}(2p)}$ can  be written as follows:
$${\rm With\ } z_{j,i} := -z_{i,j} {\rm \ for\ } i \leq j, \ \sum_{k=1}^{2p} { \overline{z_{j,k}} }z_{i,k} = \delta_{i,j},\ \ 1 \leq i \leq j \leq 2p$$
Formally replacing $z_{i,j}$ by $S_{i,j}$ and ${\overline {z_{i,j}}}$ by $S_{i,j}^*$ (where $S_{i,j} \in {\mathcal B}({\mathcal H}))$, one is led to formulate Theorems 5.1 and 5.2 below.\\

{\bf Theorem 5.1}. Let $S= (S_{1,1}, \ldots,S_{1,p};S_{2,2},\ldots,S_{2,p};\ldots; S_{p,p})$ be a $\frac{p(p+1)}{2}$-tuple of commuting operators in ${\mathcal B}({\mathcal H})$. Then (a) and (b) below are equivalent.
\\
(a) $S$ is an $S_{\Omega_{II}(p)}$-isometry.\\
(b) With $S_{j,i} := S_{i,j}$ for $i \leq j$,
$$\sum_{k=1}^p S_{j,k}^*S_{i,k} = \delta_{i,j}I_{\mathcal H},\ \ 1 \leq i \leq j \leq p.$$\\
\begin{proof} The necessity of the conditions (b) is by now obvious. For the sufficiency part we note that the conditions in (b) corresponding to $1 \leq i=j \leq p$ guarantee that each $S_{l,m}$, with $l \leq m$, is an operator coordinate of a $p$-tuple that is a spherical isometry so that Lemma 3.4 applies. One can then argue as in the proof of Theorem 4.1.
\end{proof}

{\bf Theorem 5.2}. Let $S= (S_{1,2}, \ldots,S_{1,2p};S_{2,3};\ldots,S_{2,2p};\ldots;S_{2p-1,2p})$ be a 
$p(2p-1)$-tuple of commuting operators in ${\mathcal B}({\mathcal H})$. Then (a) and (b) below are equivalent.
\\
(a) $S$ is an $S_{\Omega_{III}(2p)}$-isometry.\\
(b) With $S_{j,i} := -S_{i,j}$ for $i \leq j$,
$$\sum_{k=1}^{2p} S_{j,k}^*S_{i,k} = \delta_{i,j}I_{\mathcal H},\ \ 1 \leq i \leq j \leq 2p.$$\\
\begin{proof} The necessity of the conditions (b) is obvious. For the sufficiency part we note that the conditions in (b) corresponding to $1 \leq i=j \leq 2p$ guarantee that each $S_{l,m}$, with $l < m$, is an operator coordinate of a $(2p-1)$-tuple that is a spherical isometry so that Lemma 3.4 applies. One can then argue as in the proof of Theorem 4.1.
\end{proof}

{\bf Remark 5.3}. In view of Theorems 3.2, 4.1, 5.1 and 5.2, it is clear that the argument in the proof of Theorem 3.5 can be pushed through to accommodate the domains $\Omega_{II}(p)$ and $\Omega_{III}(2p)$ as well and  the statement of Theorem 4.2 stands generalized by way of letting each $\Omega_i$ to be any of $\Omega_I(p,q)$, $\Omega_{IV}(n)$, $\Omega_{II}(p)$ and $\Omega_{III}(2p)$.\\

We now turn our attention to the domains $\Omega_{III}(2p+1)$. The Shilov boundary $S_{\Omega_{III}(2p+1)}$ is the set
$$ \{z_{{\mathcal A}(2p+1)} \in \C^{p(2p+1)}: {\rm With\ } z_{j,i} := -z_{i,j} {\rm \  for \ } i \leq j ,\  (z_{i,j}) = UKU^t {\rm \  for \ some\ } {\rm \ unitary \ matrix\ } U\},$$
where
$$ K = \underbrace{\left[{\begin{array}{cc} 0 & 1\\ -1 & 0 \end{array}} \right] \oplus \cdots \oplus \left[{\begin{array}{cc} 0 & 1\\ -1 & 0 \end{array}} \right]}_{\hbox{$p$ summands}} \oplus[0].$$\\
The matrix $Z := (z_{i,j})= UKU^t$ is such that $Z^*Z$ has $0$ as a characteristic value of multiplicity $1$ and $1$ as a characteristic value of multiplicity $2p$.\\

For $p(2p+1)$-tuples $z_{{\mathcal A}(2p+1)}$ and $w_{{\mathcal A}(2p+1)}$, we let $z_{j,i} = -z_{i,j}$, $w_{j,i} = -w_{i,j}$ for $i \leq j$ and, for the $(2p+1)\times(2p+1)$ antisymmetric 
matrices $Z=(z_{i,j})$ and $W=(w_{i,j})$, we let $q(\lambda; Z, W)$ denote the characteristic
polynomial $\det(\lambda I_{2p+1}-W^tZ)$ of $W^tZ$. We write $q(\lambda; Z, W)$ as
$$ q(\lambda; Z, W) = q_0(Z,W)+q_1(Z,W)\lambda +\cdots + q_{2p+1}(Z,W)\lambda^{2p+1}.$$
Any $q_k(Z,W)$ is a polynomial in the $2p(2p+1$) variables $z_{1,2},\cdots,z_{2p,2p+1},w_{1,2},\cdots,w_{2p,2p+1}$.\\

{\bf Theorem 5.4}. Let $S= (S_{1,2}, \ldots,S_{1,2p+1};S_{2,3},\ldots,S_{2,2p+1};\ldots; S_{2p,2p+1})$ be a $p(2p+1)$-tuple of commuting operators in ${\mathcal B}({\mathcal H})$. Then (a) and (b) below are equivalent.\\
(a) $S$ is an $S_{\Omega_{III}(2p+1)}$-isometry.\\
(b) $$q_0(Z,W)(S,S^*) = 0_{\mathcal H};\ q_m(Z,W)(S,S^*) = (-1)^{m-1}{{2p} \choose {m-1}}I_{\mathcal H},\ 1 \leq m \leq 2p+1.$$
\begin{proof} Suppose $S$ is an $S_{\Omega_{III}(2p+1)}$-isometry. Then the Taylor spectrum
$\sigma(N)$ of the minimal normal extension $N= (N_{1,2}, \ldots,N_{1,2p+1};N_{2,3},\ldots,N_{2,2p+1};\ldots; N_{2p,2p+1}) \in {\mathcal B}({\mathcal K})^{p(2p+1)}$ of $S$ is contained in $S_{\Omega_{III}(2p+1)}$. Since for any 
$z_{{\mathcal A}(2p+1)} \in S_{\Omega_{III}(2p+1)}$ the matrix  $Z^*Z$ has  $0$ as a characteristic value of multiplicity $1$ and $1$ as a characteristic value of multiplicity $2p$, the characteristic polynomial $q(\lambda; Z, \bar Z)$ of $Z^*Z$ coincides with $\lambda(\lambda-1)^{2p}$ and the scalar equalities
$$q_0(Z,\bar Z) = 0;\ q_m(Z,\bar Z) = (-1)^{m-1}{{2p} \choose {m-1}},\ 1 \leq m \leq 2p+1$$
hold. The operator equalities
$$q_0(Z,W)(N,N^*) = 0_{\mathcal K};\ q_m(Z,W)(N,N^*) = (-1)^{m-1}{{2p} \choose {m-1}}I_{\mathcal K},\ 1 \leq m \leq 2p+1$$
follow. Compressing the last equations to $\mathcal H$, (b) is seen to hold.

Conversely, suppose (b) holds. The condition $q_{2p}(Z,W)(S,S^*) = -2pI_{\mathcal H}$ gives 
$$ S_{1,2}^*S_{1,2}+\cdots+ S_{2p,2p+1}^*S_{2p,2p+1} = pI_{\mathcal H}$$
so that $(1/\sqrt p)S$ is a spherical isometry. It follows that $(1\sqrt p)S$ and hence $S$ is subnormal. Let $N$ in the notation used above be the minimal normal extension of $S$. Now Lemma 3.3 yields 
$$q_0(Z,W)(N,N^*) = 0_{\mathcal K};\ q_m(Z,W)(N,N^*) = (-1)^{m-1}{{2p} \choose {m-1}}I_{\mathcal K},\ 1 \leq m \leq 2p+1.$$ 
By the spectral theory for $N$, the scalar equalities
$$q_0(Z,\bar Z) = 0;\ q_m(Z,\bar Z) = (-1)^{m-1}{{2p} \choose {m-1}},\ 1 \leq m \leq 2p+1$$
hold for any $z_{{\mathcal A}(2p+1)}$ in the Taylor spectrum $\sigma(N)$ of $N$. But then the characteristic polynomial $q(\lambda; Z, \bar Z)$ of $Z^*Z$ coincides with $\lambda(\lambda-1)^{2p}$ so that $Z^*Z$ has  $0$ as a characteristic value of multiplicity $1$ and $1$ as a characteristic value of multiplicity $2p$. At this stage, we invoke a result originally due to Hua \cite{Hu0} (see also \cite[THEOREM 1]{SW}) to assert the existence of a unitary matrix $U$ such that $UZU^t = K$. But this clearly implies  $z_{{\mathcal A}(2p+1)} \in S_{\Omega_{III}(2p+1)}$.
\end{proof}

{\bf  Remark 5.5}. As observed in the proof of Theorem 5.5, any $S_{\Omega_{III}(2p+1)}$-isometry $S$ is such that $(1/\sqrt p)S$ is a spherical isometry. This necessitates, for our purposes, that the following elementary observation be made: Suppose $S_i$ is an $n_i$-tuple of operators in ${\mathcal B}({\mathcal H})$ for $1 \leq i \leq m$ with $S=(S_1;\ldots;S_m)$ being an $(n_1+\cdots+n_m)$-tuple of commuting operators. If the set $\{1,\ldots,m\}$ can be partitioned into sets $\{p_1,\ldots,p_k\}$ and $\{q_1,\ldots,q_l\}$ such that each $S_{p_i}$ satisfies the hypotheses of Lemma 3.4 and each $S_{q_j}$ is such that $(1/m_j)S_{q_j}$ is a spherical isometry for some positive number $m_j$, then $S$ is subnormal. Indeed, the tuple $S'$ consisting of $S_{p_i}$ and $(1/m_j)S_{q_j}$ satisfies the hypotheses of Lemma 3.4 and hence admits a normal extension $N$ with commuting coordinates $N_{p_i}$ and $N_{q_j}$; the tuple $N$ with the coordinates $N_{p_i}$ and $m_jN_{q_j}$ is then a normal extension of $S$.\\

Using Theorems 3.2, 4.1, 5.1, 5.2, 5.4, Remark 5.5 and arguing as in the proof of Theorem 3.5, one can now establish Theorem 5.6 below.\\

{\bf Theorem 5.6}.  Let $\Omega = \Omega_1\times\cdots\times\Omega_m \subset \C^{n}$ where each $\Omega_i$ is a classical Cartan domain of any of the types I, II, III and IV in $\C^{n_i}$ (and where $n =n_1+\cdots+n_m$). Let $S_i=(S_{i,1},\ldots,S_{i,n_i})$ be an $n_i$-tuple of operators in ${\mathcal B}({\mathcal H})$ for $1 \leq i \leq m$ and let the operator coordinates of the $n$-tuple $S=(S_1;\ldots;S_m)$ commute with each other. Then $S$ is an $S_{\Omega}$-isometry if and only if each $S_i$ is an $S_{\Omega_i}$-isometry.\\

It is interesting to note how the ``stars-on-the-left" functional calculus, in conjunction with the known characterization of an $S_{\B_n}$-isometry as a spherical isometry, 
facilitates our arguments in Sections 3, 4 and 5.

\end{document}